\documentclass[preprint,12pt]{elsarticle}

\usepackage{amssymb}

\usepackage{amsmath}
\usepackage{mathtools}
\usepackage{subcaption}
\usepackage{hyperref}

\usepackage{lineno}

\usepackage{fancyhdr}

\newcommand\numberthis{\addtocounter{equation}{1}\tag{\theequation}}

\makeatletter
\DeclareRobustCommand{\pder}[1]{%
  \@ifnextchar\bgroup{\@pder{#1}}{\@pder{}{#1}}}
\newcommand{\@pder}[2]{\frac{\partial#1}{\partial#2}}
\makeatother

\makeatletter
\DeclareRobustCommand{\der}[1]{%
  \@ifnextchar\bgroup{\@der{#1}}{\@der{}{#1}}}
\newcommand{\@der}[2]{\frac{d#1}{d#2}}
\makeatother

\newcommand{\RR}{\mathbb{R}} 

\newcommand{\overbar}[1]{\mkern 1.5mu\overline{\mkern-1.5mu#1\mkern-1.5mu}\mkern 1.5mu}

\newcommand*\diff{\mathop{}\!\mathrm{d}}

\usepackage{xcolor}
\newcommand{\rev}[1]{{#1}}
\newcommand{\myrev}[1]{{#1}}

\journal{Journal of Computational Physics}
\fancypagestyle{titlepagestyle}
{
   \fancyhf{}
   \fancyfoot[L]{\footnotesize\itshape Preprint submitted to Journal of Computational Physics\\
   \textcopyright 2019. This manuscript version is made available under the CC-BY-NC-ND 4.0 license \url{http://creativecommons.org/licenses/by-nc-nd/4.0/}}
   \fancyfoot[R]{\footnotesize\itshape June 17, 2019}
   
}

\begin{document}

\begin{frontmatter}



\title{A block preconditioner for non-isothermal flow in porous media}


\author[1]{Thomas Roy\corref{cor1}}
\ead{thomas.roy@maths.ox.ac.uk}
\cortext[cor1]{Corresponding author}

\author[2]{Tom B. J\"{o}nsth\"{o}vel}
\ead{tjonsthovel@slb.com}

\author[3]{Christopher Lemon}
\ead{clemon@slb.com}

\author[1]{Andrew J. Wathen}
\ead{andy.wathen@maths.ox.ac.uk}

\address[1]{Mathematical Institute, University of Oxford, Oxford, United Kingdom}
\address[2]{Schlumberger Houston Production Technology Center, Houston, TX, United States of America}
\address[3]{Schlumberger Abingdon Technology Center, Abingdon, United Kingdom}

\begin{abstract}
 In petroleum reservoir simulation, the industry standard preconditioner, the constrained pressure residual method (CPR), is a two-stage process which involves solving a restricted pressure system with Algebraic Multigrid (AMG). Initially designed for isothermal models, this approach is often used in the thermal case. However, it does not have a specific treatment of the additional energy conservation equation and temperature variable. We seek to develop preconditioners which better capture thermal effects such as heat diffusion. In order to study the effects of both pressure and temperature on fluid and heat flow, we consider a model of non-isothermal single phase flow through porous media. For this model, we develop a block preconditioner with an efficient Schur complement approximation. Both the pressure block and the approximate Schur complement are approximately inverted using an AMG V-cycle. The resulting solver is scalable with respect to problem size and parallelization.
\end{abstract}

\begin{keyword}
preconditioning  \sep iterative solvers \sep porous media \sep thermal reservoir simulation
\end{keyword}

\end{frontmatter}
\thispagestyle{titlepagestyle}

\section{Introduction}\label{sec:intro}
Models of fluid flow in porous media are used in the simulation of applications such as petroleum reservoirs, carbon storage, hydrogeology, and geothermal energy. In some cases, fluid flow must be coupled with heat flow in order to capture thermal effects. Petroleum reservoir simulation is used in optimizing oil recovery processes, which often involve heating and steam injection inside the reservoir in order to reduce the viscosity of the oil. This is especially important in the case of heavier hydrocarbons. 

In the case of isothermal multiphase flow, a global pressure couples local concentration/saturation variables. The equations in the system are elliptic with respect to the pressure and hyperbolic with respect to the non-pressure variables. The industry standard constrained pressure residual (CPR) preconditioner introduced by Wallis ~\cite{wallis1983incomplete,wallis1985constrained} in the early 80s defines a discrete decoupling operator essentially splitting pressure and non-pressure variables, in order that each can be preconditioned separately. Indeed, the global nature of the pressure variable requires a more precise ``global'' preconditioning than the other variables for which ``local'' preconditioning is sufficient. In brief, the CPR preconditioner is a two-stage process in which pressure is solved first approximately, followed by solving approximately the full system. 

A major improvement to CPR was introduced in \cite{lacroix2003iterative} with the use of Algebraic Multigrid (AMG) \cite{ruge1987algebraic} as a preconditioner for the pressure equation in the first stage. AMG is used as a solver for elliptic problems, usually as a preconditioner for a Krylov subspace method. Therefore, the elliptic-like nature of the pressure equation makes it an ideal candidate for the use of AMG. This improved preconditioner, often denoted CPR-AMG, is widely used in modern reservoir simulators. 

The non-isothermal case adds a conservation of energy equation and a temperature (or enthalpy) variable to the system of PDEs. In the standard preconditioning approach, the energy conservation equation and the temperature unknowns are treated similarly to the secondary equations and unknowns. This means that the thermal effects are only treated in the second stage of CPR, usually an Incomplete LU factorization (ILU) method. More dense incomplete factors are often needed in the thermal case. While this results in a lower iteration count, it is not ideal in terms of computational time, memory requirements, and parallelization.  

Alternatives to the usual approach were recently proposed. The Fraunhofer Institute for Algorithms and Scientific Computing (SCAI) focuses on AMG for systems of PDEs based on \cite{clees2005amg}, often called System AMG (SAMG). While it seeks to replace CPR-AMG, the SAMG approach proposed in \cite{clees2010efficient,gries2014preconditioning} is quite similar in the isothermal case. Indeed, AMG is applied to the whole system, but all non-pressure variables remain on the fine level. In the thermal case, however, SAMG allows the consideration of both pressure and temperature for the hierarchy. A proper comparison with CPR-AMG for thermal simulation cases has yet to be done. \rev{Other AMG methods for systems of PDEs include BoomerAMG \cite{henson2002boomeramg} and multigrid reduction (MGR) in the hypre library \cite{falgout2002hypre}, as well as Smoothed Aggregation in the ML package \cite{ml-guide}. In particular, BoomerAMG has been shown to be effective in diffusion-dominated two-phase flow problems \cite{bui2017algebraic}, and MGR has also had some success with multiphase flow problems \cite{wang2017multigrid, bui2018algebraic}. }

The inclusion of temperature in the AMG hierarchy is still not well understood. The temperature is not always descriptive of the flow everywhere in the reservoir (for example in regions of faster flows). This justifies an adaptive method where only variables which are descriptive be included in the first stage of CPR. Retaining the CPR structure, Enhanced CPR (ECPR) constructs a ``strong'' subsystem for the first stage of CPR by looking at the coupling in the system matrix \cite{li2017enhanced}. The resulting subsystem has no real physical interpretation and it is unclear if it is possible to solve it via AMG. 

In the context of this paper, we consider single phase non-isothermal flow. This single phase case is relevant for geothermal models and simple reservoir simulation examples, but can also be applied to miscible displacement problems (where a concentration plays a similar role to temperature) \cite{booth2008miscible}. We present a block preconditioner for the solution of the resulting coupled pressure-temperature system. 

In Section \ref{sec:problem}, we present the mathematical model for non-isothermal flow in porous media and the discretization. In Section \ref{sec:precons}, we describe the preconditioning approaches for the linearized system. Numerical results for the preconditioners are presented in Section \ref{sec:results}. We conclude in Section \ref{sec:conclusion} with a discussion on the future direction of this research. 


\section{Problem statement}
\label{sec:problem}
In this section, we describe a coupled PDE system and its discretization. 
\subsection{Single phase thermal flow in porous media}
\label{sec:thermaleqns}

We describe the equations for single phase flow in porous media coupled with thermal effects.

\subsubsection{Conservation of mass}

We start with the continuity equation which states that the rate at which mass enters the system is equal to the rate of mass which leaves the system plus the accumulation of mass within the system. Additionally, we include a source/sink term which accounts for mass which is added or removed from the system. We have
\begin{equation}\label{eq:continuity1}
 \phi\pder{\rho}{t} + \nabla \cdot (\rho \mathbf{u}) = f \quad \text{in } \RR_+\times\Omega,
\end{equation}
where $\phi$ is the porosity field of the rock, $\rho$ is the density of the fluid, $\mathbf{u}$ is the fluid velocity, $f$ is a source/sink term, and $\Omega$ is the spatial domain in $\RR^d$, $d=2,3$. The source/sink term $f$ represents injection/production wells and is given in Section \ref{sec:source}. We further assume that the velocity follows Darcy's law \cite{darcy1856fontaines}, i.e.
\begin{equation}\label{eq:darcy}
  \mathbf{u} = -\frac{\mathbf{K}}{\mu}(\nabla p - \rho \mathbf{g}),
\end{equation}
where $p$ is the pressure, $\mathbf{K}$ is the permeability tensor field, $\mu$ is the viscosity, and $\mathbf{g}$ is gravitational acceleration. The density and viscosity are functions of pressure and temperature given in Section \ref{sec:murho}. Then, \eqref{eq:continuity1} becomes
\begin{equation}\label{eq:consmass}
 \phi\pder{\rho}{t} - \nabla \cdot \left(\rho \frac{\mathbf{K}}{\mu}(\nabla p - \rho \mathbf{g})\right) = f \quad \text{in } \RR_+\times\Omega.
\end{equation}
We also assume Neumann and Dirichlet boundary conditions
 \begin{equation}\label{eq:massbc}
  -\frac{\mathbf{K}}{\mu}(\nabla p - \rho \mathbf{g})\cdot \mathbf{n} = g_N \text{ on } \Gamma_N,\quad\text{and}\quad p = g_D  \text{ on } \Gamma_D,
 \end{equation}
where $g_N$ is Neumann boundary data, $g_D$ is Dirichlet boundary data, $\mathbf{n}$ is the unit outward normal vector on $\partial \Omega = \Gamma_N \cup \Gamma_D$, and $\Gamma_D \cap \Gamma_N = \varnothing$. 

\subsubsection{Conservation of energy}

Similarly, we have a conservation of energy equation for the heat energy. Note that formulations where enthalpy is an independent variable are common, but here we consider temperature as an independent variable as in a reference commercial reservoir simulator \cite{debaun2005extensible}.
Here, $c_v$ and $c_r$ are the specific heat of the fluid and rock, respectively, $\rho_r$ is the density of the rock, and $T$ is temperature. Here, $c_v T$ represents the enthalpy of the fluid, and $\rho c_v $, its energy density. Heat energy is not only transported by a heat flux, but also by the fluid flux. We get the following advection-diffusion equation:
\begin{equation}\label{eq:energy1}
 \phi\pder{}{t}(\rho c_v T) + (1-\phi)\pder{}{t}(\rho_r c_r T) + \nabla \cdot (\rho c_v T \mathbf{u}) + \nabla \cdot\mathbf{q} = f_T \quad  \text{in } \RR_+\times\Omega,
\end{equation}
where $\mathbf{q}$ is the heat flux, and $f_T$ is a source/sink term representing wells or heaters and given in Section \ref{sec:source}. Furthermore, we assume that the heat flux follows Fourier's law, i.e.
\begin{equation}
 \mathbf{q} = -k_T \nabla T,
\end{equation}
where $k_T$ is the thermal conductivity field. It is given by 
\begin{equation}
 k_T = \phi k_{T,r} + (1-\phi) k_{T,f},
\end{equation}
where $k_{T,r}$ and $k_{T,f}$ are the conductivities of the rock and the fluid, respectively. Then, \eqref{eq:energy1} becomes
\begin{equation}\label{eq:energy2}
 \phi\pder{}{t}(\rho c_v T) + (1-\phi)\pder{}{t}(\rho_r c_r T) + \nabla \cdot (\rho c_v T \mathbf{u}) - \nabla \cdot (k_T \nabla T) = f_T \quad  \text{in } \RR_+\times\Omega.
\end{equation}
Then, assuming Darcy flow, we get
\begin{multline}\label{eq:energy3}
 \phi\pder{}{t}(\rho c_v T) + (1-\phi)\pder{}{t}(\rho_r c_r T) -\nabla \cdot \left(\rho c_v T \frac{\mathbf{K}}{\mu}(\nabla p - \rho \mathbf{g})\right) - \nabla \cdot (k_T \nabla T) = f_T \\  \text{in } \RR_+\times\Omega.
\end{multline}
We also assume Neumann and Dirichlet boundary conditions
\begin{equation}\label{eq:energy4}
  -\left( \rho c_v T \frac{\mathbf{K}}{\mu}(\nabla p - \rho \mathbf{g}) + k_T\nabla T \right) \cdot \mathbf{n} = g^T_N \text{ on } \Gamma^T_N,\quad\text{and}\quad T = g^T_D \text{ on } \Gamma^T_D,
\end{equation}
where $g_N^T$ is Neumann boundary data, $g_D^T$ is Dirichlet boundary data, $\partial \Omega = \Gamma^T_N \cup \Gamma^T_D$, and $\Gamma^T_D \cap \Gamma^T_N = \varnothing$.


\subsubsection{Coupled problem}
We assume that $\rho$ and $\mu$ are empirically determined functions of pressure and temperature. Our choices are given in Section \ref{sec:murho}.

We are interested in solving the following boundary value problem:
\\
find $p$, $T$ such that
\begin{equation}\label{eq:strong1}
 \phi\pder{\rho}{t} - \nabla \cdot \left(\rho \frac{\mathbf{K}}{\mu}(\nabla p - \rho \mathbf{g})\right) = f \quad \text{in } \RR_+\times\Omega,
 \end{equation}
 \begin{multline}\label{eq:strong2}
   \phi\pder{}{t}(\rho c_v T) + (1-\phi)\pder{}{t}(\rho_r c_r T) -\nabla \cdot \left(\rho c_v T \frac{\mathbf{K}}{\mu}(\nabla p - \rho \mathbf{g})\right) - \nabla \cdot (k_T \nabla T) = f_T \\ \text{in } \RR_+\times\Omega,
 \end{multline}
 \begin{equation}\label{eq:strong3}
  -\frac{\mathbf{K}}{\mu}(\nabla p - \rho \mathbf{g})\cdot \mathbf{n} = g_N \text{ on } \Gamma_N,\quad\text{and}\quad p = g_D  \text{ on } \Gamma_D,
 \end{equation}
\begin{equation}\label{eq:strong4}
 -\left( \rho c_v T \frac{\mathbf{K}}{\mu}(\nabla p - \rho \mathbf{g}) + k_T\nabla T \right) \cdot \mathbf{n} = g^T_N \text{ on } \Gamma^T_N,\quad\text{and}\quad T = g^T_D \text{ on } \Gamma^T_D,
\end{equation}
where $\partial \Omega = \Gamma_N \cup \Gamma_D =  \Gamma_N^T \cup \Gamma_D^T$, $\Gamma_D \cap \Gamma_N =\Gamma_D^T \cap \Gamma_N^T = \varnothing$, and initial conditions for $p$ and $T$ are prescribed. 

\subsubsection{Nonlinear quantities}\label{sec:murho}
The density $\rho$ and viscosity $\mu$ are empirically determined functions of temperature and pressure. For the examples in this paper, we will consider the flow of heavy oil in porous media and thus use the following empirical laws.

\begin{table}[hbt]
\centering
\caption{Parameters for the Bennison viscosity correlation}\label{tab:bennison}
 \begin{tabular}{c c c c}
 \hline
  $A_1$ & $A_2$ & $A_3$ & $A_4$ \\
  \hline
  -0.8021 & 23.8765 & 0.31458 & -9.21592\\
  \hline
 \end{tabular}
\end{table}

\begin{figure}[hbt]
 \centering
 \includegraphics[width = 0.9\textwidth]{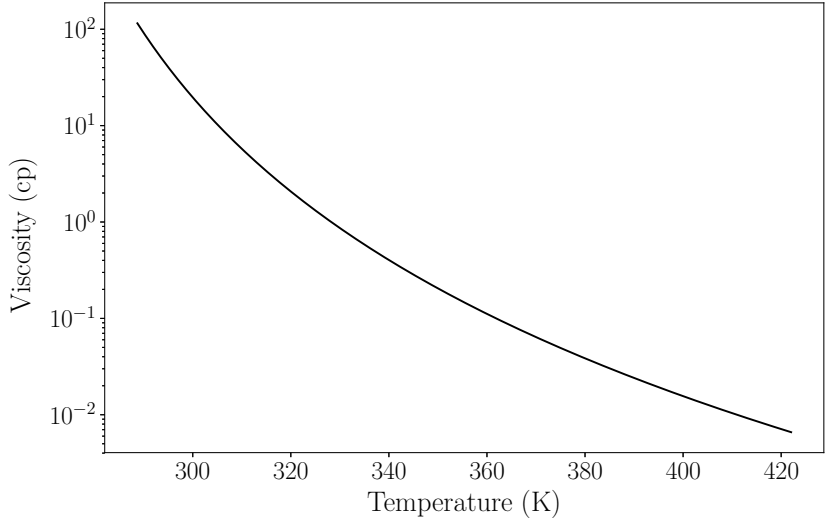}
 \caption{The Bennison viscosity correlation for heavy oil.}\label{fig:bennison}
\end{figure}

For viscosity, we choose the following correlation \cite{bennison1998prediction}:
\begin{equation}\label{eq:bennison}
 \mu (T_\mathrm F) = 10^{A_1 \gamma_\mathrm{API} + A_2} T_\mathrm F^{A_3 \gamma_\mathrm{API} + A_4},
\end{equation}
which takes temperature $T_\mathrm F$ in $^\circ$F and returns viscosity in cp (0.001 kg m$^{-1}$ s$^{-1}$). The viscosity as a function of temperature (in Kelvin) is illustrated in Figure \ref{fig:bennison}. The dimensionless parameters $A_i$ can be found in Table \ref{tab:bennison}. The American Petroleum Institute (API) gravity $\gamma_\mathrm{API}$ is a measure of how heavy or light a petroleum liquid is compared to water: if its API is greater than 10, than it is lighter and floats on water; if less than 10, it is heavier and sinks. We can calculate API gravity from specific gravity (SG) (ratio of the density of the petroleum liquid to the density of water, at 60$^\circ$ F) using the following formula:
\begin{equation}
 \gamma_\mathrm{API} = \frac{141.5}{\gamma_\mathrm{SG}} - 131.5.
\end{equation}
For density, we use the following correlation:
\begin{equation}
 \rho(p,T) = \rho_0 e^{c (p-p_0)} e^{\beta (T - T_0)},
\end{equation}
where $p_0$, $T_0$ are reference pressure and temperature and $\rho_0$ is the density at those values, $c$ is a compressibility coefficient and $\beta$ is a thermal expansion coefficient. 
 Values representative to those used in reservoir simulation are $p_0=1.01325$ bar, $T_0$= 288.7056 K ($60 \;^\circ$F), $c= 5.5\times 10^{-5} \;\mathrm{bar}^{-1}$, and $\beta = 2.5\times 10^{-4}\; \mathrm{K}^{-1}$. 
Given a specific gravity, we have $\rho_0 = \gamma_\mathrm{SG}\; \rho_w$, where $\rho_w = 999$ kg m$^{-3}$ is the density of water at the reference temperature.

\subsubsection{Source/sink terms}\label{sec:source}
We first consider source/sink terms representing injection and production wells. A simple way to model these is by using point sources/sinks
\begin{align*}
 f(\mathbf x ) =& \sum_i q^i_\mathrm{inj}(p,T) \delta (\mathbf x - \mathbf x_\mathrm{inj}^i) \rho  (p,T_\mathrm{inj}) \\
	     &-\sum_j q^j_\mathrm{prod}(p,T) \delta (\mathbf x - \mathbf x_\mathrm{prod}^j) \rho (p,T) , \numberthis \label{eq:f}
\end{align*}
\begin{align*}
 f_T(\mathbf x ) =& \sum_i q^i_\mathrm{inj}(p,T) \delta (\mathbf x - \mathbf x_\mathrm{inj}^i) \rho  (p,T_\mathrm{inj}) c_v T_\mathrm{inj} \\
	     &-\sum_j q^j_\mathrm{prod} (p,T)\delta (\mathbf x - \mathbf x_\mathrm{prod}^j) \rho (p,T) c_v T  , \numberthis \label{eq:fT}
\end{align*}
where $\mathbf x_\mathrm{inj}$ and $\mathbf x_\mathrm{prod}$ represent the location of injection and production wells, respectively, $\delta (\mathbf{x})$ is the Dirac delta function, $q_\mathrm{inj}^i$ and $q_\mathrm{prod}^j$ are the wells' injection and production rates, respectively. 

The production rate $q_\mathrm{prod}$ is usually given by a constant target production rate. Similarly, the injection rate $q_\mathrm{inj}$ is given by a target injection rate. These rates can only be maintained if the pressure at the production well does not drop below a minimum pressure, and the pressure at the injection well does not go above a maximum pressure. In those cases, a well model is required. We consider the commonly used Peaceman well model \cite{peaceman1978interpretation,chen2009well} for anisotropic media with $\mathbf{K} = \mathrm{diag} (K_x, K_y, K_z)$ as the permeability tensor field. In this case the rates are given by
\begin{equation}
 q = \frac{2\pi h K_e}{\mu \mathrm{ln} (r_e/r_w)} (p_{bh} - p),
\end{equation}
where $h$ is the height of well opening, $K_e= \sqrt{K_x K_y}$ is the equivalent permeability, $p_{bh}$ is the bottom-hole pressure, $r_w$ is the well radius, and $r_e$ is the equivalent radius which can be calculated using
\begin{equation}
 r_e = \frac{0.14\left( (K_y/K_x)^{1/2} D_x^2 + (K_x/K_y)^{1/2} D_y^2)\right)^{1/2}}{0.5\left((K_y/K_x)^{1/4} + (K_x/K_y)^{1/4}\right)},
\end{equation}
where $D_x$ and $D_y$ are the horizontal lengths of the grid cell. Since we want to allow mesh refinements, we do not want the model to change as we vary the grid size. Therefore, we arbitrarily fix $D_x = D_y = 5$ meters, and also choose $h=5$ meters and $r_w=0.1$ meters.

Oil recovery techniques for heavy oils can include electromagnetic heating \cite{sahni2000electromagnetic}. These can be expressed as source terms for the energy equation. For simplicity, we do not use an electromagnetic model and choose the simple function 
\begin{equation}
 f_T =  \sum_i U_\mathrm{heater}(p,T) \delta (\mathbf x - \mathbf x_\mathrm{heater}^i) (T_\mathrm{heater} - T),
\end{equation}
where $\mathbf x_\mathrm{heater}$ represent the location of heaters, $U_\mathrm{heater}$ is the heat transfer coefficient, and $T_\mathrm{heater}$ is the target heating temperature. For our simulations, we have a heating coefficient of $5.44409\times 10 ^{-6}$ Js$^{-1}$K$^{-1}$. For simplicity, we also choose $T_\mathrm{heater}$ to be the same as $T_\mathrm{inj}$. 


\subsection{DG0 discretization}
\label{sec:fem}

In reservoir simulation, Finite Volume methods are most commonly used \cite{leveque2002finite}. Since the flux entering a given volume is identical to that leaving an adjacent one, these methods are conservative. Additionally, upwind schemes introduce substantial numerical diffusion, which helps with stability. In this section, we present a discontinuous Galerkin (DG) method \cite{riviere2008discontinuous} that is equivalent to a Finite Volume method used in reservoir simulation and is based on the description in \cite{riseth2015nonlinear}. The resulting weak formulation allows us to implement our problem in the open source Finite Element software Firedrake \cite{rathgeber2016firedrake}.

Let $\mathcal{T} = \{E_i, i\in \mathcal{I}\}$ be a partition of $\Omega$ into open element domains $E_i$ such that union of their closure is $\overbar \Omega$, where $\mathcal{I}$ is a set of indices. Let the interior facet $e_{ij} = \overbar E_i\cap \overbar E_j$ and let $\Gamma_\mathrm{int}$ denote the union of all interior facets. Let connection set $\mathcal{N}(i)$ denote the set of indices $j$ such that $|e_{ij}|>0$. We begin by presenting a DG0 (piecewise constant) method for the heat equation 
\begin{equation}\label{eq:heat}
 \pder{u}{t} -\nabla^2 u = 0 \text{ in } \Omega,
\end{equation}
\begin{equation}\label{eq:heat2}
 u = f  \text{ on } \Gamma_D,\quad \nabla u \cdot \mathbf n = g \text{ on } \Gamma_N.
\end{equation}
The variational problem for \eqref{eq:heat}-\eqref{eq:heat2} on a single cell $E_i$ is: find $u$ such that
\begin{equation}\label{eq:heatweak}
\int_{E_i} \pder{u}{t} v \diff x + \int_{E_i} \nabla u \cdot \nabla v \diff x - \int_{\partial E_i} v \nabla u \cdot \mathbf n \diff s = 0 \quad \text{for all test functions } v,
\end{equation}
where $\mathbf n$ is the outward normal to $E_i$. Let us first consider $E_i$ such that $\partial E_i \in \Gamma_\mathrm{int}$. Let $h_i$ denote the center point of cell $E_i$. For the flux on the interior facets, we choose the following flux approximation
\begin{equation}
 \int_{\partial E_i} v \nabla u \cdot \mathbf n \diff s \approx \sum_{j\in \mathcal{N} (i)} \int_{e_{ij}} v\mid_{E_i}^{e_{ij}} \frac{u\mid_{E_j}^{e_{ij}} - u\mid_{E_i}^{e_{ij}}}{\|h_j - h_i\|}\diff s.
\end{equation}
Here, $u\mid_{E_i}^{e_{ij}}$ denotes the limit of $u$ in cell $E_i$ as it goes to the edge $e_{ij}$. 

We consider a piecewise constant approximation of our solution, i.e. in the approximation space $\mathcal{V}_h = \mathbb P_\mathrm{DG}^0$ with basis $\{\phi_i = \mathbf{1}_{E_i}\mid i\in \mathcal{I} \}$. The DG0 approximation is $u_h = \sum_{i\in\mathcal{I}} u_i \phi_i$. For this approximation, on $E_i$, $v_h\in \mathcal V_h$ is constant and $\nabla v_h = 0$. Therefore, \eqref{eq:heatweak} becomes
\begin{equation}
 \int_{E_i} \pder{u_i}{t} v_i \diff x - \sum_{j\in \mathcal{N} (i)}  \int_{e_{ij}} v_i \frac{u_j - u_i}{\|h_j - h_i\|}\diff s = 0.
\end{equation}
Note that this is equivalent to 
\begin{equation}
 \pder{u_i}{t} |E_i|- \sum_{j\in \mathcal{N} (i)}  \frac{u_j - u_i}{\|h_j - h_i\|}  |e_{ij}| = 0,
\end{equation}
which is a Finite Volume approximation of the heat equation. In reservoir simulation, this way of approximating the interior facet integrals is known as a ``two-point flux'' (TPFA) approximation. In order for such a Finite Volume method to converge, the grid must satisfy a certain orthogonality property \cite{eymard2000finite}. In brief, in each cell, there exists a point called the center of the cell such that for any adjacent cell, the straight line between the two centers is orthogonal to the boundary between the cells. For the examples in this paper, we choose quadrilateral meshes, which easily satisfy this condition.

If instead $E_i$ is a boundary element, then the boundary integral becomes
\begin{equation}
 \int_{\partial E_i} v \nabla u \cdot \mathbf n \diff s = \sum_{j\in \mathcal{N} (i)} \int_{e_{ij}} v \nabla u \cdot\mathbf n_e  + \int_{\partial E_i \cap \Gamma_D} v \nabla u \cdot\mathbf n \diff s  + \int_{\partial E_i \cap \Gamma_N} v \nabla u \cdot\mathbf n \diff s,
\end{equation}
where $\mathbf n_e$ is the unit outward pointing normal of a cell. We use the following flux approximation
\begin{align*}
 \int_{\partial E_i} v \nabla u \cdot \mathbf n \diff s \coloneqq& \sum_{j\in \mathcal{N} (i)} \int_{e_{ij}}  v_i \frac{u_j - u_i}{\|h_j - h_i\|}\diff s \\ 
 &+ \int_{\partial E_i \cap \Gamma_D}  v_i \frac{(f - u_i)}{d_{\Gamma_D}(h_i)}\diff s + \int_{\partial E_i \cap \Gamma_N} v_i g\diff s, \numberthis
\end{align*}
where $d_{\Gamma_D}(h_i)$ is the shortest distance form $h_i$ to the boundary $\Gamma_D$. For each $i\in \mathcal{I}$, we have
\begin{align*}
 \int_{E_i} \pder{u_i}{t} v_i \diff x - \sum_{j\in \mathcal{N} (i)}  \int_{e_{ij}} v_i \frac{u_j - u_i}{\|h_j - h_i\|}\diff s -\int_{\partial E_i \cap \Gamma_D}  v_i \frac{(f - u_i)}{d_{\Gamma_D}(h_i)}\diff s 
 \\- \int_{\partial E_i \cap \Gamma_N} v_i g\diff s = 0 . \numberthis
\end{align*}
For a given ordering of the indices in $\mathcal{I}$, we denote by $u^+$ and $u^-$ the limit value of $u$ for two cells sharing an edge. Now, summing over all $i\in\mathcal{I}$, and noting that each interior facet is visited twice, we obtain
\begin{equation}
 \int_\Omega \pder{u}{t} v\diff x + \int_{\Gamma_\mathrm{int}} (v^+ - v^-) \frac{u^+ - u^-}{\|h^+ - h^-\|} \diff S - \int_{\Gamma_D} v \frac{f-u}{d_{\Gamma_D}(h)} \diff s- \int_{\Gamma_N} v g \diff s = 0.
\end{equation}
We define the jump of $v$ as $[v] = v^+ - v^-$. We then get the following problem: find $u\in \mathbb P_\mathrm{DG}^0$ such that 
\begin{equation}
 \int_\Omega \pder{u}{t} v\diff x + \int_{\Gamma_\mathrm{int}} [v]\frac{[u]}{\|h^+ - h^-\|} \diff S - \int_{\Gamma_D} v \frac{f-u}{d_{\Gamma_D}(h)} \diff s- \int_{\Gamma_N} v g \diff s = 0,
\end{equation}
for all $v\in \mathbb P_\mathrm{DG}^0$.

\subsubsection{Upwinding}

We now consider an upwind Godunov method \cite{godunov1959difference} for the advection equation
\begin{equation}
 \pder{u}{t} + \nabla \cdot (u \mathbf{w}) = 0 \text{ on } \Omega,
\end{equation}
\begin{equation}
 u = f  \text{ on } \Gamma_D,\quad  u \mathbf w\cdot \mathbf n = g \text{ on } \Gamma_N,
\end{equation}
where $\mathbf{w}$ is a given vector field. For an interior $E_i$, the upwind scheme is given by
\begin{equation}
 \int_{E_i} \pder{u}{t} v \diff x + \int_{\partial E_i} v u^\mathrm{up} \mathbf{w} \cdot \mathbf n \diff s = 0,
\end{equation}
where $u^\mathrm{up}$ is the upwind value of $u$, which, for a facet $e$ shared by $E_1$ and $E_2$ and $\mathbf n_e$ pointing from $E_1$ to $E_2$, is given by
\begin{equation}
 u^\mathrm{up} = \begin{cases}
                  u\mid_{E_1}^e	&\text{if } \mathbf w \cdot \mathbf n_e \geq 0,\\
                  u\mid_{E_2}^e	&\text{if } \mathbf w \cdot \mathbf n_e < 0.
                 \end{cases}
\end{equation}
For the full discretized problem we have: find $u\in \mathbb P_\mathrm{DG}^0$ such that 
\begin{equation}
 \int_\Omega \pder{u}{t} v \diff x + \int_{\Gamma_\mathrm{int}} [v] u^\mathrm{up} \mathbf{w} \cdot \mathbf n \diff s - \int_{\Gamma_D} v \frac{f-u}{d_{\Gamma_D}(h)} \diff s- \int_{\Gamma_N} v g \diff s = 0,
\end{equation}
for all $v\in \mathbb P_\mathrm{DG}^0$.

\subsubsection{Semidiscrete problem}
We now discretize \eqref{eq:strong1}-\eqref{eq:strong4} in space using the semidiscrete DG0 formulation described above. Assuming homogeneous Neumann boundary conditions, the variational problem is: find the approximation $(p,T)\in \mathbb P_\mathrm{DG}^0 \times \mathbb P_\mathrm{DG}^0$ such that 
\begin{align*}
\int_\Omega \myrev{\phi}\pder{\rho}{t} q \diff x +  \int_{\Gamma_\mathrm{int}} [q] \left( \{\!\!\{\mathbf{K}\}\!\!\}\frac{\rho^\mathrm{up}}{\mu^\mathrm{up}}\left(\frac{[p]}{\|h^+ - h^-\|} - \{\rho\} \mathbf{g} \cdot\mathbf n_e\right)\right)\diff S \\- \int_\Omega f q \diff x = 0, \numberthis
\end{align*}
\begin{multline}
 \int_\Omega \myrev{\phi} c_v\pder{\rho T}{t} r \diff x
 + \int_\Omega \myrev{(1-\phi)}\rho_r c_r \pder{T}{t} r \diff x \\
 +  \int_{\Gamma_\mathrm{int}} [r]  \{\!\!\{\mathbf{K}\}\!\!\}c_v\frac{\rho^\mathrm{up}}{\mu^\mathrm{up}}T^\mathrm{up}\left(\frac{[p]}{\|h^+ - h^-\|} - \{\rho\} \mathbf{g} \cdot\mathbf n_e\right)\diff S \\
 + \int_{\Gamma_\mathrm{int}} [r] \{\!\!\{k_T\}\!\!\}\frac{[T]}{\|h^+ - h^-\|}\diff S- \int_\Omega f_T r \diff x = 0,
\end{multline}
for all $(q,r)\in \mathbb P_\mathrm{DG}^0 \times \mathbb P_\mathrm{DG}^0$. The brackets $\{\}$ denote the average across the facets, and the double brackets $\{\!\!\{\}\!\!\}$ denote the harmonic average across the facets. The use of the harmonic average is standard for two-point flux approximation, and is obtained by considering piecewise constant permeabilities \cite{eymard2000finite}. The upwind quantities are given by
\begin{equation}
 (u)^\mathrm{up} = \begin{cases}
                  u\mid_{E_1}^e	&\text{if  }  \dfrac{[p]}{\|h^+ - h^-\|} - \{\rho\} \mathbf{g} \cdot\mathbf n_e \geq 0,\\
                  u\mid_{E_2}^e	&\text{if  }  \dfrac{[p]}{\|h^+ - h^-\|} - \{\rho\} \mathbf{g} \cdot\mathbf n_e < 0.
                 \end{cases}
\end{equation}
For the delta functions in the source/sink terms, we choose the simple approximation:
\begin{equation}
 \delta (x) = \begin{cases}
            1/|E_i| &\text{if } x \in E_i, \\
            0 &\text{otherwise}.
            \end{cases}
\end{equation}

\subsubsection{Fully discretized problem}
For time discretization, we use the backward Euler method. We define the two following forms, which are linear with respect with their last argument:
\begin{multline}\label{eq:Fm}
F_m(p^{n+1},T^{n+1};q) \coloneqq \int_\Omega \myrev{\phi}\frac{\rho^{n+1} - \rho^n}{\Delta t} q \diff x \\
+  \int_{\Gamma_\mathrm{int}} [q] \left( \{\!\!\{\mathbf{K}\}\!\!\}\frac{(\rho^{n+1})^\mathrm{up}}{(\mu^{n+1})^\mathrm{up}}\left(\frac{[p^{n+1}]}{\|h^+ - h^-\|} - \{\rho^{n+1}\} \mathbf{g} \cdot\mathbf n_e\right)\right)\diff S \\
- \int_\Omega f^{n+1} q \diff x,
\end{multline}
\begin{multline}\label{eq:Fe}
 F_e(p^{n+1},T^{n+1};r) \coloneqq \int_\Omega \myrev{\phi}c_v\frac{\rho^{n+1} T^{n+1} - \rho^n T^n}{\Delta t} r \diff x \\
 + \int_\Omega \myrev{(1-\phi)}\rho_r c_r \frac{T^{n+1} - T^n}{\Delta t} r \diff x\\
 +  \int_{\Gamma_\mathrm{int}} [r]  \{\!\!\{\mathbf{K}\}\!\!\}\frac{(\rho^{n+1})^\mathrm{up}}{(\mu^{n+1})^\mathrm{up}}(T^{n+1})^\mathrm{up}\left(\frac{[p^{n+1}]}{\|h^+ - h^-\|} - \{\rho^{n+1}\} \mathbf{g} \cdot\mathbf n_e\right)\diff S \\
 + \int_{\Gamma_\mathrm{int}} [r] \{\!\!\{k_T\}\!\!\}\frac{[T^{n+1}]}{\|h^+ - h^-\|}\diff S- \int_\Omega f_T^{n+1} r \diff x.
\end{multline}
Let $ F(p,T;q,r) \coloneqq F_m(p,T;q) +  F_e(p,T;r)$, which is linear in both $q$ and $r$, but nonlinear in $p$ and $T$. At each time-step, given the previous solution $(p^n,T^n) $, we search for $(p^{n+1},T^{n+1})\in\mathbb P_\mathrm{DG}^0\times\mathbb P_\mathrm{DG}^0$ such that
\begin{equation}\label{eq:nonlinear}
 F(p^{n+1},T^{n+1};q,r) = 0\quad \text{for all } (q,r)\in \mathbb P_\mathrm{DG}^0\times\mathbb P_\mathrm{DG}^0.
\end{equation}


\section{Solution algorithms}
\label{sec:precons}

The system of nonlinear equations \eqref{eq:nonlinear} can be written as a system of nonlinear equations for the real coefficients $p_i$ and $T_i$ of the DG0 functions $p^{n+1}$ and $T^{n+1}$, respectively. Let $x$ be the vector of these coefficients and $G$ the function such that $G(x) = 0$ is equivalent to \eqref{eq:nonlinear}. By linearizing this equation with Newton's method, we must solve at each iteration
\begin{equation}
 \pder{G}{x}\mid_{x=x_k} (x_{k+1} - x_k) = -G(x_k).
\end{equation}
The resulting linearized systems can be written as a block system of the form
\begin{equation}\label{eq:AXR}
A \delta x = \begin{bmatrix}
A_{pp} & A_{pT} \\
A_{Tp} & A_{TT}
\end{bmatrix}
\begin{bmatrix}
\delta p \\ \delta T
\end{bmatrix}
= \begin{bmatrix}
b_p \\ b_T
\end{bmatrix} = b,
\end{equation}
where $\delta x= x_{k+1}-x_k$ is the Newton increment. The different blocks are the discrete versions of Jacobian terms as follows

\begin{equation}\label{eq:App}
 A_{pp} \sim \phi \frac{1}{\Delta t} (\rho)_p + \nabla \cdot \left(\rho \mathbf{u}\right)_p - (f)_p,
\end{equation}
\begin{equation}
 A_{pT} \sim \phi \frac{1}{\Delta t} (\rho)_T  + \nabla \cdot \left(\rho \mathbf{u}\right)_T - (f)_T,
\end{equation}
\begin{equation}
 A_{Tp} \sim \phi \frac{1}{\Delta t} (\rho)_p c_v T + \nabla \cdot \left(c_v T(\rho \mathbf{u})_p\right) - (f_T)_p,
\end{equation}
\begin{multline}
  A_{TT} \sim \phi \frac{c_v (\rho + (\rho)_T T) }{\Delta t} + (1-\phi) \frac{\rho_r c_r }{\Delta t} + \nabla \cdot \left( c_v \rho \mathbf{u}\right) \\ + \nabla \cdot \left( c_v  T (\rho \mathbf{u})_T\right) - \nabla \cdot (k_T \nabla) - (f_T)_T,
\end{multline}
where 
\begin{equation}
 (\rho \mathbf{u})_p = -\frac{\mathbf K}{\mu} \left(\rho (\nabla - (\rho)_p \mathbf{g})+ (\rho)_p (\nabla p - \rho \mathbf{g} )\right),
\end{equation}
and
\begin{equation}\label{eq:rhouT}
 (\rho \mathbf{u})_T = -\mathbf K \left[\left( \frac{\rho}{\mu} \right)_T (\nabla p - \rho \mathbf{g} ) \myrev{- \frac{\rho}{\mu} (\rho)_T \mathbf{g} }\right].
\end{equation}
All coefficients in \eqref{eq:App}-\eqref{eq:rhouT} are evaluated at the previous Newton iterate $(p_k, T_k)$, and $(.)_p$ and $(.)_T$ denote the partial derivatives with respect to $p$ and $T$, respectively.

The linearized systems are often very difficult to solve using iterative methods. Indeed, efficient preconditioning is required in order to achieve rapid convergence with linear solvers \cite{saad2003iterative}. In this section, we will detail different preconditioning techniques used to solve \eqref{eq:AXR}. We first mention some methods which are important ingredients of the preconditioning techniques. 

Krylov subspace methods are used to approximate the solution of $Ax = b$ by constructing a sequence of Krylov subspaces, $\mathcal{K}_n = \left\{b,Ab, A^2 b, \dots, A^{n-1} b\right\}$. The generalized minimal residual method (GMRES) \cite{saad1986gmres} is a Krylov subspace method suitable for general linear systems. The approximate solution $x_n$ is formed by minimizing the Euclidean norm of the residual $r_n = Ax_n -b$ over the subspace $\mathcal{K}_n$.

Incomplete LU factorization (ILU) \cite{saad2003iterative,meijerink1977iterative} is a general preconditioning technique in which sparse triangular factors are used to approximate the system matrix $A$. This preconditioner requires assembling the factors and then solving two triangular systems. A popular way to determine the sparsity pattern of the factors is to simply choose the relevant triangular parts of the sparsity pattern of $A$. This is known as ILU(0). More generally, choosing the sparsity pattern of $A^{k+1}$ is called ILU($k$).

Multigrid methods \cite{brandt1977multi,briggs2000multigrid, trottenberg2000multigrid} use hierarchies of coarse grid approximations in order to solve differential equations. Smoothing operations (such as a Jacobi or Gauss-Seidel iteration) are combined with coarse grid corrections on increasingly coarser grids. For positive definite elliptic PDEs, it is known that multigrid methods can provide optimal solvers (in the sense of linear scalability with the dimension of the discretized problem). 

Algebraic Multigrid (AMG) \cite{ruge1987algebraic,stuben2001introduction} uses information from the entries of the system matrix rather than that of the geometric grid. This makes AMG an ideal black-box solver for elliptic problems. Although it can be used to solve simpler problems, AMG is often used as a preconditioner for Krylov subspace methods in problems which are essentially elliptic. Relative to preconditioners such as ILU, parallel variants of multigrid methods retain more effectiveness.


\subsection{Two-stage preconditioning: CPR}
\label{sec:twostage}

Let $M_1$ and $M_2$ be two preconditioners for the linear system $Ax = b$ for which we have the action of their (generally approximate) inverse $M_1^{-1}$, and $M_2^{-1}$. Applying a multiplicative two-stage preconditioner can be done as follows:
\begin{enumerate}
 \item Precondition using $M_1$: $x_1 = M_1^{-1}b$;
 \item Compute the new residual: $b_1 = b - A x_1$;
 \item Precondition using $M_2$ and correct: $x = M_2^{-1} b_1 + x_1$.
\end{enumerate}
The action of the two-stage preconditioner can be written as
\begin{equation}
 M^{-1} = M^{-1}_{2} (I - AM_{1}^{-1}) + M^{-1}_1.
\end{equation}
In the case of multiphase flow in porous media, the standard preconditioner is the Constrained Pressure Residual method (CPR) \cite{wallis1983incomplete}. In the multiphase case, the linear systems are like in \eqref{eq:AXR} except that the temperature blocks are replaced (or combined in the thermal case) with saturation blocks. 

In the case of CPR, the first stage preconditioner $M_1$ is given by
\begin{equation} 
    M_1^{-1} \approx \begin{bmatrix}
                            A_{pp}^{-1} & 0 \\
                            0 & 0
                            \end{bmatrix},
\end{equation}
 where $A_{pp}^{-1}$ is approximated using an AMG V-cycle. The second preconditioner is chosen such that $M_2^{-1} \approx A^{-1}$, usually with an incomplete LU factorization method (ILU).

 In addition to the two-stage preconditioner, decoupling operators are often used to reduce the coupling between the pressure equation and the saturation variables. Indeed, an approximation of the pressure equation $A_{pp} \delta p + A_{pT} \delta T= b_T$ is performed in the first stage of CPR where the saturation coupling $A_{pT}$ is ignored. A decoupling operator is a left preconditioner applied a priori to \eqref{eq:AXR} of the form
\begin{equation}
 \begin{bmatrix}
       I & -D \\
       0 & I 
      \end{bmatrix}.
\end{equation}
The most often-used approximations for multiphase flow are Quasi-IMPES (QI) and True-IMPES (TI) ~\cite{lacroix2003iterative,lacroix2000iterative,scheichl2003decoupling}. The approximations are $D_{QI} =\mathrm{diag} (A_{pT}) \mathrm{diag} (A_{TT})^{-1}$, $D_{TI} =\mathrm{colsum} (A_{pT}) \mathrm{colsum} (A_{TT})^{-1}$. Here, $\mathrm{colsum}(A)$ is a diagonal matrix with entries the sums of the entries in the columns of $A$, which is equivalent to the mass accumulation terms when discretized with the two-point flux approximation as outlined in this paper (as the fluxes sum up to zero in a given column). 

By performing this decoupling operation on the system \eqref{eq:AXR} before CPR, the first stage now consists in solving a subsystem for the approximate Schur complement $S_p = A_{pp} - D A_{Tp}$ instead of the original pressure block. However, the properties of the resulting $S_p$ need to be amenable to the application of AMG (for example M-matrix properties). While this is nearly guaranteed in the black-oil case \cite{gries2015system}, it does not necessarily follow for compositional flow or thermal flow. For the single phase test cases detailed in Section \ref{sec:results}, we observe that CPR performs best without decoupling operators (results not shown here).

\subsection{Block factorization preconditioner}
\label{sec:block}
Consider the following decomposition of the Jacobian 

 \begin{equation}\label{eq:Afact}
 A = 
  \begin{bmatrix}
   I & 0 \\
   A_{Tp} A_{pp}^{-1} & I
  \end{bmatrix}
  \begin{bmatrix}
   A_{pp} & 0 \\
   0 & S_T
  \end{bmatrix}
  \begin{bmatrix}
   I & A_{pp}^{-1} A_{pT} \\
   0 & I
  \end{bmatrix},
 \end{equation}
 where $S_T = A_{TT} - A_{Tp} A_{pp}^{-1} A_{pT}$ is the Schur complement. The inverse of the Jacobian is given by
 
  \begin{equation}\label{eq:Ainv}
 A^{-1} = 
  \begin{bmatrix}
   I & -A_{pp}^{-1} A_{pT} \\
   0 & I
  \end{bmatrix}
  \begin{bmatrix}
   A_{pp}^{-1} & 0 \\
   0 & S_T^{-1}
  \end{bmatrix}
  \begin{bmatrix}
   I & 0 \\
   - A_{Tp} A_{pp}^{-1} & I
  \end{bmatrix}.
 \end{equation}
 Even if $A$ is sparse, the Schur complement $S_T$ is generally dense. A common preconditioning technique is to use the blocks of the factorization \eqref{eq:Afact} combined with a sparse approximation of the Schur complement \cite{elman2014finite}. \rev{Given an appropriate Schur complement approximation $\tilde{S}_T$, applying the block preconditioner can be done as follows:
 \begin{enumerate}
  \item Solve the pressure subsystem: $A_{pp} x_p = b_p$;
  \item Compute the new energy equation residual: $\tilde{b}_T = b_T - A_{Tp} x_p$;
  \item Solve the Schur complement subsystem: $\tilde{S}_T \delta_T = \tilde{b}_T$;
  \item Compute the new mass equation residual: $\tilde{b}_p = x_p - A_{pT} \delta_T$;
  \item Solve the pressure subsystem: $A_{pp} \delta_p = \tilde{b}_p$.
 \end{enumerate}
In our case, $A_{pp}^{-1}$ and $\tilde{S}_T^{-1}$ are both approximated using an AMG V-cycle. }

\subsection{Schur complement approximation} 
 Common sparse approximations for the Schur complement are $\myrev{\tilde{S}_{A_{TT}}} = A_{TT}$ and $\tilde{S}_{\mathrm{diag}} = A_{TT} - A_{Tp} \mathrm{diag}\left( A_{pp}\right)^{-1} A_{pT}$. Here we present a Schur complement approximation which performs significantly better than such simple approximations. 
 
 For the derivation of our Schur complement approximation, we consider the linearized problem before discretization. This approach results in an approximation which holds as we refine the mesh. See \cite{mardal2011preconditioning} for a theoretical framework in using the infinite-dimensional setting to find mesh-independent preconditioners for self-adjoint problems. 

\subsubsection{Steady-state case} %
\label{sec:schur}
We first consider a steady-state single phase thermal problem: find $p$, $T$ such that
\begin{equation}\label{eq:steady1}
\nabla \cdot \left(\rho \mathbf u \right) = 0 \quad \text{in } \Omega,
 \end{equation}
 \begin{equation}\label{eq:steady2}
 \nabla \cdot \left(\rho c_v T \mathbf u \right) - \nabla \cdot (k_T \nabla T) = 0 \quad  \text{in } \Omega,
 \end{equation}
where $\mathbf{u}$ is given by \eqref{eq:darcy}, and we have no-flux boundary conditions for the fluid and heat. Here we will consider the linearized system in a continuous setting. Applying a Newton method to \eqref{eq:steady1}-\eqref{eq:steady2}, we obtain a block systems of the form \eqref{eq:AXR} where the blocks are:
\begin{equation}
 A_{pp} = \nabla \cdot (\rho \mathbf u)_p, \qquad A_{pT} = \nabla \cdot (\rho \mathbf u)_T ,
\end{equation}
\begin{equation}\label{eq:atp}
 A_{Tp} = \nabla \cdot \left(c_v T (\rho \mathbf u)_p \right) = c_v \left[\nabla T \cdot (\rho \mathbf u)_p + T \nabla \cdot (\rho \mathbf u)_p \right], 
\end{equation}
 \begin{align*}
 A_{TT} &= c_v \left[ \nabla \cdot (\rho \mathbf u) + \nabla \cdot \left( T (\rho \mathbf u)_T \right) \right] - \nabla \cdot (k_T \nabla) \\
 &= c_v \left[ \nabla \cdot (\rho \mathbf u) + \nabla T \cdot (\rho \mathbf u)_T + T \nabla \cdot  (\rho \mathbf u)_T  \right] - \nabla \cdot (k_T \nabla), \numberthis \label{eq:att}
\end{align*}
where we have used the product rule for the divergence operator in \eqref{eq:atp} and \eqref{eq:att}. Then the second term of the Schur complement (which corresponds in the continuous setting to the Poincar\'e-Steklov operator) becomes
\begin{align*}
 A_{Tp} A_{pp}^{-1} A_{pT}  
 &=  c_v \left[\nabla T \cdot (\rho \mathbf u)_p + T \nabla \cdot (\rho \mathbf u)_p \right] \left(  \nabla \cdot (\rho \mathbf u)_p  \right)^{-1} \nabla \cdot (\rho \mathbf u)_T \\
 &= c_v T \nabla \cdot (\rho \mathbf u)_T + c_v\nabla T \cdot (\rho \mathbf u)_p \left(  \nabla \cdot (\rho \mathbf u)_p  \right)^{-1} \nabla \cdot (\rho \mathbf u)_T. \numberthis \label{eq:schur2}
\end{align*}
We notice that in $A_{TT} - A_{Tp} A_{pp}^{-1} A_{pT}$, the terms $c_v T \nabla \cdot (\rho \mathbf u)_T$ cancel. We are left with
\begin{multline}\label{eq:fullschur}
 S_T = c_v \nabla \cdot (\rho \mathbf u) + c_v \nabla T \cdot (\rho \mathbf u)_T - \nabla \cdot (k_T \nabla) + c_v \nabla T \cdot (\rho \mathbf u)_T \\ - c_v\nabla T \cdot (\rho \mathbf u)_p \left(  \nabla \cdot (\rho \mathbf u)_p  \right)^{-1} \nabla \cdot (\rho \mathbf u)_T .
\end{multline}
One of the nonlinear terms has canceled, and so we consider if it is possible that the last two terms also cancel. Consider the operator $(\rho \mathbf u)_p \left(  \nabla \cdot (\rho \mathbf u)_p  \right)^{-1} \nabla \cdot$, which is close to the operator $\nabla (\nabla \cdot \nabla)^{-1} \nabla \cdot \eqqcolon s$. This holds if $(\rho \mathbf u)_p$ is close to $\nabla$, i.e. if $\rho$ is close to being constant with respect to $p$. \rev{While this approximation holds for liquid water and hydrocarbons, it may be less applicable in the case of gases. Extending this to multiphase flow is straightforward and is part of ongoing work.}

Assuming that the operator $s$ is applied to a sufficiently smooth vector field $\mathbf F$, we can use Helmholtz decomposition to decompose this field into the sum of its curl-free and divergence-free part $\mathbf F = -\nabla \Phi + \nabla \times \mathbf A$, where $\Phi$ is a scalar potential and $\mathbf A$ a vector potential. Since $s$ removes the divergence-free part of a field, applying $s$ to $\mathbf F$, we obtain
\begin{equation}
 s \mathbf F = - \nabla (\nabla \cdot \nabla)^{-1} \nabla \cdot \nabla \Phi = -\nabla \Phi,
\end{equation}
assuming that $\Phi$ satisfies the same boundary conditions as the operator $(\nabla \cdot \nabla)^{-1}$. Hence $s$ is a projection to the curl-free subspace, and it acts like the identity operator when applied to curl-free vector fields.

 We assume that this also holds for $(\rho \mathbf u)_p \left(  \nabla \cdot (\rho \mathbf u)_p  \right)^{-1} \nabla \cdot$. In \eqref{eq:fullschur}, we see that this operator is applied to $(\rho \mathbf u)_T$ which is of the form $\gamma \nabla p$, where $\gamma$ is a scalar field. In order for $\gamma \nabla p$ to be a curl-free vector field, we need $\nabla \times (\gamma \nabla p) = \nabla \gamma \times \nabla p = 0$, i.e. we need $\nabla \gamma$ and $\nabla p$ to be parallel vectors. In the discretized case, our grid satisfies an orthogonality property as mentioned in Section \ref{sec:fem}, and the gradients of $p$ and $\gamma$ are approximated using a two-point flux approximation. In this case, the gradients are always orthogonal to the facets, and thus parallel. Accordingly, we replace the operator $(\rho \mathbf u)_p \left(  \nabla \cdot (\rho \mathbf u)_p  \right)^{-1} \nabla \cdot$ by the identity and obtain the following Schur complement approximation
\begin{equation}\label{eq:schurapprox}
 \tilde{S}_T = c_v \nabla \cdot (\rho \mathbf u) - \nabla \cdot (k_T \nabla).
\end{equation}
Similar heuristic arguments for replacing $\nabla \cdot (\nabla \cdot \nabla)^{-1} \nabla$ by the identity operator in the case of the Stokes problem can be found, for example, in \cite{maday1993analysis}, and for the Navier-Stokes equations, in \cite{kay2002preconditioner}.

 \subsubsection{Source terms}\label{sec:sourceschur}

Similarly, we consider the steady-state case with the addition of source/sink terms. In this case, production wells satisfy $f^\mathrm{prod}_T = c_v T f_\mathrm{prod}$, while injection wells satisfy $f^\mathrm{inj}_T = c_v T_\mathrm{inj} f_\mathrm{inj}$. Thus,

\begin{align*}
  S_{T} =& \nabla \cdot \left( c_v \rho \mathbf{u}\right)- \nabla \cdot (k_T \nabla) - c_v f_\mathrm{prod} +  c_v \nabla T  \cdot  (\rho \mathbf{u})_T \\
  &+ c_v T\nabla \cdot (\rho \mathbf{u})_T -c_v T (f_\mathrm{prod})_T - c_v T_\mathrm{inj} (f_\mathrm{inj})_T 
  \\
  & - \left(c_v\nabla T \cdot (\rho \mathbf u)_p +c_v T\nabla \cdot (\rho \mathbf u )_p -c_v T (f_\mathrm{prod})_p - c_v T_\mathrm{inj} (f_\mathrm{inj})_p \right) \\ 
  &\left( \nabla \cdot (\rho \mathbf u)_p   - (f_\mathrm{prod})_p  - (f_\mathrm{inj})_p \right)^{-1}  \left[\nabla \cdot (\rho \mathbf u)_T    - (f_\mathrm{prod})_T  - (f_\mathrm{inj})_T \right].  \numberthis\label{eq:wellschur}
\end{align*}
Since the injection term is weighted by $T_\mathrm{inj}$, we cannot directly cancel the $c_v T$ terms as in \eqref{eq:fullschur}. However, after a certain amount of injection, $T$ tends to $T_\mathrm{inj}$ where the injection well is located. Furthermore, in the infinite-dimensional setting, this effect will be instantaneous since the well terms are defined using a delta function in \eqref{eq:fT}. Using this argument, we get

\begin{align*}
  S_{T} \approx& \nabla \cdot \left( c_v \rho \mathbf{u}\right)- \nabla \cdot (k_T \nabla) - c_v f_\mathrm{prod} +  c_v \nabla T  \cdot  (\rho \mathbf{u})_T \\
  & - \left(c_v\nabla T \cdot (\rho \mathbf u)_p  \right) \\ 
  &\left( \nabla \cdot (\rho \mathbf u)_p   - (f_\mathrm{prod})_p  - (f_\mathrm{inj})_p \right)^{-1}  \left[\nabla \cdot (\rho \mathbf u)_T    - (f_\mathrm{prod})_T  - (f_\mathrm{inj})_T \right].  \numberthis \label{eq:wellschur2}
\end{align*}
Further assuming that the mass source/sink terms are almost constant in $p$ and $T$, i.e. $(\rho)_T$ and $(\rho)_p$ are small and the injection/production rates are independent of pressure and temperature (which is the case when operating at a target rate), we ignore the derivatives of the source/sink terms. Then, using the same argument as for the steady-state case, we obtain the Schur complement approximation
\begin{equation}\label{eq:wellschurapprox}
  \tilde S_{T} =\nabla \cdot \left( c_v \rho \mathbf{u}\right)- \nabla \cdot (k_T \nabla) - c_v f_\mathrm{prod} 
\end{equation}

In the case where the the source terms are heaters, we have $f=0$, and  $f_T = U (T_\mathrm{heater} - T) D_\mathrm{heaters}$, where $D_\mathrm{heaters}$ is the sum of delta functions for the location of heaters. The Schur complement is given by
\begin{align*}
  S_{T} =& \nabla \cdot \left( c_v \rho \mathbf{u}\right) - \nabla \cdot (k_T \nabla) + UD_\mathrm{heaters} +  c_v \nabla T  \cdot  (\rho \mathbf{u})_T 
  \\
  &- c_v\nabla T \cdot (\rho \mathbf u)_p \left( \nabla \cdot (\rho \mathbf u)_p  \right)^{-1}  \left[\nabla \cdot (\rho \mathbf u)_T\right].
\end{align*}
We see that heaters do not affect the right-hand side term. Using the same argument as above, we get the approximation
\begin{equation}
  \tilde S_{T} = \nabla \cdot \left( c_v \rho \mathbf{u}\right)- \nabla \cdot (k_T \nabla) + U D_\mathrm{heaters}.
\end{equation}
 
\subsubsection{Time-dependent case}

We now generalize our analysis to the time-dependent problem. We first consider the case without source/sink terms. The blocks are given by 

\begin{equation}
 A_{pp} = \phi \frac{1}{\Delta t} (\rho)_p + \nabla \cdot \left(\rho \mathbf{u}\right)_p,
\end{equation}
\begin{equation}
 A_{pT} = \phi \frac{1}{\Delta t} (\rho)_T  + \nabla \cdot \left(\rho \mathbf{u}\right)_T,
\end{equation}
\begin{equation}
 A_{Tp} = \phi \frac{1}{\Delta t} (\rho)_p c_v T + \nabla \cdot \left(c_v T(\rho \mathbf{u})_p\right),
\end{equation}
\begin{multline}
  A_{TT} = \phi \frac{c_v (\rho + (\rho)_T T) }{\Delta t} + (1-\phi) \frac{\rho_r c_r }{\Delta t} + \nabla \cdot \left( c_v \rho \mathbf{u}\right) \\ + \nabla \cdot \left( c_v  T (\rho \mathbf{u})_T\right) - \nabla \cdot (k_T \nabla).
\end{multline}
The second term of the Schur complement is given by
\begin{align*}
 A_{Tp} A_{pp}^{-1} A_{pT}  
 =&  c_v \left[ \phi \frac{1}{\Delta t} (\rho)_p c_v T+ \nabla T \cdot (\rho \mathbf u)_p + T \nabla \cdot (\rho \mathbf u)_p \right] 
 \\ &\left( \phi \frac{1}{\Delta t} (\rho)_p +  \nabla \cdot (\rho \mathbf u)_p   \right)^{-1} 
 \left[ \phi \frac{1}{\Delta t} (\rho)_T  + \nabla \cdot (\rho \mathbf u)_T  \right]\\
 =& c_v T \left[ \nabla \cdot (\rho \mathbf u)_T + \phi \frac{1}{\Delta t} (\rho)_T \right] + c_v\nabla T \cdot (\rho \mathbf u)_p \\
 &\left( \phi \frac{1}{\Delta t} (\rho)_p +  \nabla \cdot (\rho \mathbf u)_p   \right)^{-1}  \left[ \phi \frac{1}{\Delta t} (\rho)_T  + \nabla \cdot (\rho \mathbf u)_T  \right], \numberthis \label{eq:schurtransient}
\end{align*}
and thus the Schur complement is
\begin{align*}
  S_{T} =& \phi \frac{c_v \rho}{\Delta t} + (1-\phi) \frac{\rho_r c_r }{\Delta t} + \nabla \cdot \left( c_v \rho \mathbf{u}\right)  - \nabla \cdot (k_T \nabla) +  c_v \nabla T  \cdot  (\rho \mathbf{u})_T 
  \\
  &- c_v\nabla T \cdot (\rho \mathbf u)_p \left( \phi \frac{1}{\Delta t} (\rho)_p +  \nabla \cdot (\rho \mathbf u)_p   \right)^{-1}  \left[ \phi \frac{1}{\Delta t} (\rho)_T  + \nabla \cdot (\rho \mathbf u)_T  \right].\numberthis
\end{align*}
To justify further simplification, we need to assume that either $\rho$ is almost constant in $p$ and $T$, or that the time-step is very large. We get the following Schur complement approximation:
\begin{equation}
  \tilde S_{T} = \phi \frac{c_v \rho}{\Delta t} + (1-\phi) \frac{\rho_r c_r }{\Delta t} + \nabla \cdot \left( c_v \rho \mathbf{u}\right)  - \nabla \cdot (k_T \nabla) +  c_v \nabla T  \cdot  (\rho \mathbf{u})_T .
\end{equation}
In the case where we have source/sink terms, the Schur complement is given by 
\begin{align*}
  S_{T} =& \phi \frac{c_v \rho}{\Delta t} + (1-\phi) \frac{\rho_r c_r }{\Delta t} + \nabla \cdot \left( c_v \rho \mathbf{u}\right)  - \nabla \cdot (k_T \nabla) +  c_v \nabla T  \cdot  (\rho \mathbf{u})_T 
  \\
  &+ c_v T \nabla \cdot (\rho \mathbf u)_T+ UD_\mathrm{heaters}  - c_v f_\mathrm{prod} -c_v T (f_\mathrm{prod})_T - c_v T_\mathrm{inj} (f_\mathrm{inj})_T 
  \\
  &- c_v \left[ \phi \frac{1}{\Delta t} (\rho)_p T+ \nabla T \cdot (\rho \mathbf u)_p + T \nabla \cdot (\rho \mathbf u)_p -T (f_\mathrm{prod})_p - T_\mathrm{inj} (f_\mathrm{inj})_p\right] 
 \\ &\left( \phi \frac{1}{\Delta t} (\rho)_p +  \nabla \cdot (\rho \mathbf u)_p   - (f_\mathrm{prod})_p  - (f_\mathrm{inj})_p \right)^{-1} \\
 &\left[ \phi \frac{1}{\Delta t} (\rho)_T  + \nabla \cdot (\rho \mathbf u)_T  - (f_\mathrm{prod})_T  - (f_\mathrm{inj})_T \right].\numberthis
\end{align*}
Using the same argument for the injection temperature $T_\mathrm{inj}$ as in Section \ref{sec:sourceschur}, we get
\begin{align*}
  S_{T} \approx& \phi \frac{c_v \rho}{\Delta t} + (1-\phi) \frac{\rho_r c_r }{\Delta t} + \nabla \cdot \left( c_v \rho \mathbf{u}\right)  - \nabla \cdot (k_T \nabla) +  c_v \nabla T  \cdot  (\rho \mathbf{u})_T 
  \\
  &+ UD_\mathrm{heaters} -c_v f_\mathrm{prod} - c_v \left[\nabla T \cdot (\rho \mathbf u)_p\right] 
  \\
  &\left( \phi \frac{1}{\Delta t} (\rho)_p +  \nabla \cdot (\rho \mathbf u)_p   - (f_\mathrm{prod})_p  - (f_\mathrm{inj})_p \right)^{-1} \\
 &\left[ \phi \frac{1}{\Delta t} (\rho)_T  + \nabla \cdot (\rho \mathbf u)_T  - (f_\mathrm{prod})_T  - (f_\mathrm{inj})_T \right].\numberthis
\end{align*}
Then, again assuming that the mass source/sink terms are independent of $p$ and $T$, we obtain the following Schur complement approximation:
\begin{equation}\label{eq:schurapproxfull}
  \tilde{S}_{T} = \phi \frac{c_v \rho}{\Delta t} + (1-\phi) \frac{\rho_r c_r }{\Delta t} + \nabla \cdot \left( c_v \rho \mathbf{u}\right)  - \nabla \cdot (k_T \nabla) + UD_\mathrm{heaters} -c_v f_\mathrm{prod}.
\end{equation}
We can obtain the discretized version of this operator from \eqref{eq:Fe} by removing the terms depending on the previous time-step, and evaluating the nonlinear terms at the previous Newton iteration. We get the following bilinear operator:
\begin{align*}
 S_e(\delta T, r) \coloneqq& \int_\Omega \myrev{\phi} c_v\frac{\rho \delta T}{\Delta t} r \diff x 
 + \int_\Omega \myrev{(1-\phi)}\rho_r c_r \frac{\delta T}{\Delta t} r \diff x\\
 &+  \int_{\Gamma_\mathrm{int}} [r]  \{\!\!\{\mathbf{K}\}\!\!\}\frac{(\rho)^\mathrm{up}}{(\mu)^\mathrm{up}}(\delta T)^\mathrm{up}\left(\frac{[p]}{\|h^+ - h^-\|} - \{\rho\} \mathbf{g} \cdot\mathbf n_e\right)\diff S \\
 &+ \int_{\Gamma_\mathrm{int}} [r] \{\!\!\{k_T\}\!\!\}\frac{[\delta T]}{\|h^+ - h^-\|}\diff S + \int_\Omega (-c_v f_\mathrm{prod} + U D_\mathrm{heaters}) \delta T\diff x. \numberthis \label{eq:Se}
\end{align*}


\section{Numerical results}
\label{sec:results}

In this section, we perform numerical experiments for our block preconditioner and CPR. These are implemented on the open source Finite Element software Firedrake \cite{rathgeber2016firedrake}. The linear algebra backend is the PETSc library \cite{petsc-web-page}, allowing efficient and parallel computations. The CPR preconditioner (without decoupling) is implemented by providing PETSc options. Our custom block preconditioner is implemented through Firedrake's Python interface. Recent work from \cite{kirby2018solver} allows us to easily assemble our Schur complement approximation preconditioner by providing a weak form with the bilinear operator \eqref{eq:Se}. We modified the custom preconditioner class from \cite{kirby2018solver} to allow the use of matrix formats other than \texttt{matfree}. For example, the default \texttt{aij} matrix format allows for faster computations for lower order methods such as the one described in Section \ref{sec:fem}. Our implementation is available on GitHub\footnote{\url{https://github.com/tlroy/thermalporous}}.

For the block preconditioner, we use our Schur complement approximation \eqref{eq:Se}, \rev{unless stated otherwise}. Both the pressure block $A_{pp}$ and the approximate Schur complement are inverted using a V-cycle of \rev{ AMG}. We use BoomerAMG \cite{henson2002boomeramg} from the hypre library \cite{falgout2002hypre} \rev{with default parameters, i.e. a symmetric-SOR/Jacobi relaxation scheme (one sweep up, one sweep down), Falgout coarsening,
classical Ruge-St\"uben interpolation, and Gaussian Elimination as the coarse grid solver.} This implementation has a very efficient parallel version of AMG. 

For the second stage of CPR, we use ILU(0) as provided from PETSc. In parallel, we use block Jacobi with ILU(0) for each block (the partition is assigned when Firedrake does the discretization).

The nonlinear solver is Newton's method with line search, and the linear solver is \rev{right-preconditioned GMRES \cite{saad1986gmres}}, restarted after 30 iterations. The convergence tolerance of Newton's methods is $10^{-8}$ for the relative function norm and relative step size norm. The convergence tolerance for GMRES is $10^{-10}$ for the relative residual norm \myrev{for the tests in Section \ref{sec:spe10}, and $10^{-5}$ for the tests in Sections \ref{sec:probsizescaling} and \ref{sec:parallel}.}

For all cases, we consider a heavy oil with density and viscosity as described in Section \ref{sec:murho}. The other physical parameters are shown in Table \ref{tab:parameters}. These parameters are representative with those used in commercial reservoir simulators. 

\begin{table}[hbt]
\centering
 \caption{Physical parameters for test cases} \label{tab:parameters}
 \begin{tabular}{l l}
 \hline
 Initial pressure       & 4.1369$\times 10^5$ Pa        \\ 
 Conductivity of oil    & 0.15 Wm$^{-1}$K$^{-1}$        \\ 
 Conductivity or rock   & 1.7295772056 Wm$^{-1}$K$^{-1}$\\ 
 Specific heat of oil   & 2093.4 JK$^{-1}$kg$^{-1}$     \\ 
 Specific heat of rock  & 920 JK$^{-1}$kg$^{-1}$        \\ 
 \hline
 \end{tabular}
\end{table}

For all cases, we evaluate the performance of the methods by comparing the number of linear iterations per nonlinear iteration. We note that, for our proof-of-concept implementation, the cost of applying the \myrev{b}lock preconditioner is around two times more computationally expensive (in serial) than CPR. The difference may not be as significant in an optimized implementation.

\subsection{SPE10 test cases}
\label{sec:spe10}
\begin{figure}[hbt]
\centering
 \begin{subfigure}[b]{0.49\textwidth}
 \includegraphics[width=\linewidth]{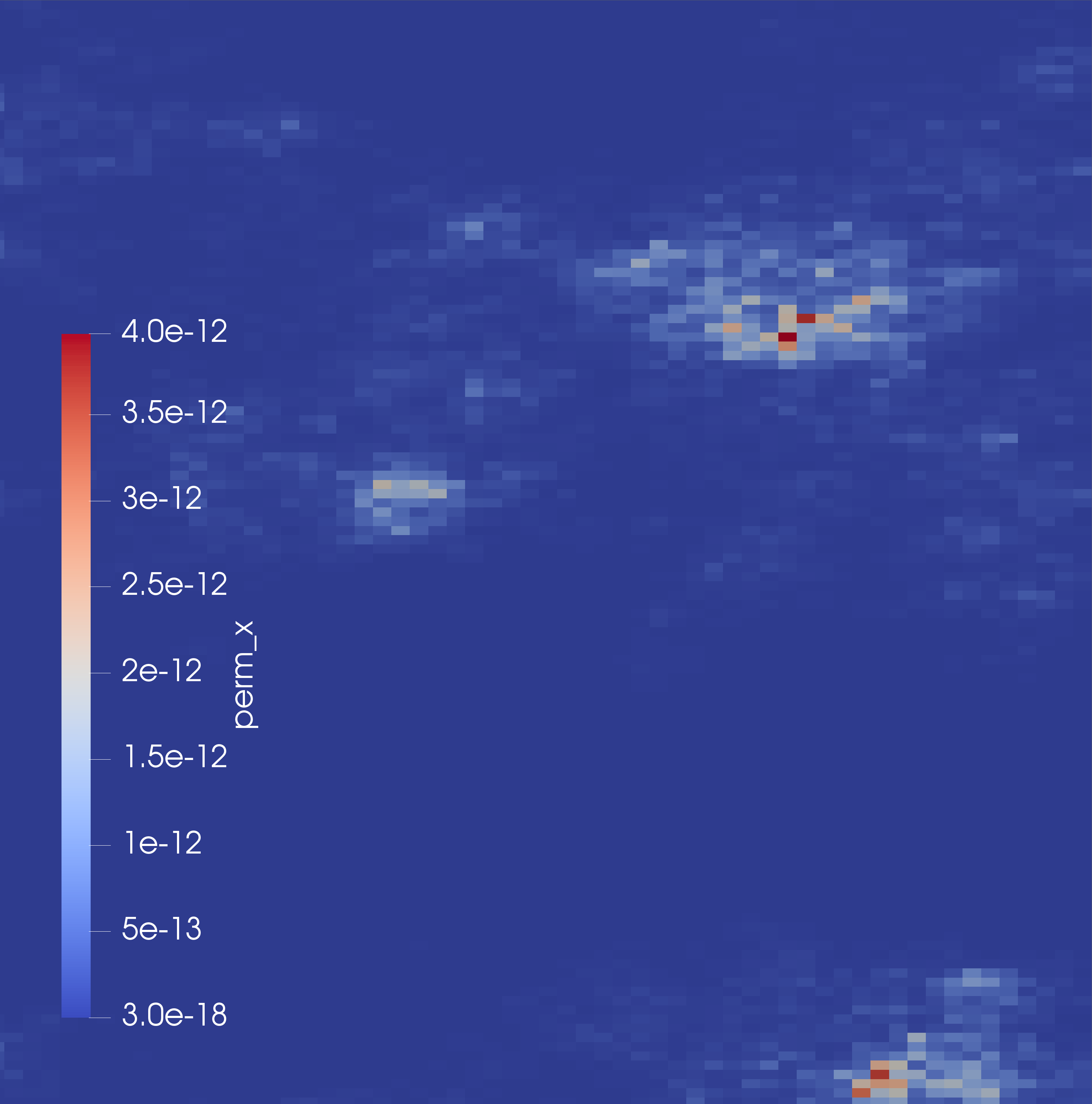}
 \caption{Linear scale}
 \end{subfigure}
 \hfill
 \begin{subfigure}[b]{0.49\textwidth}
 \includegraphics[width=\textwidth]{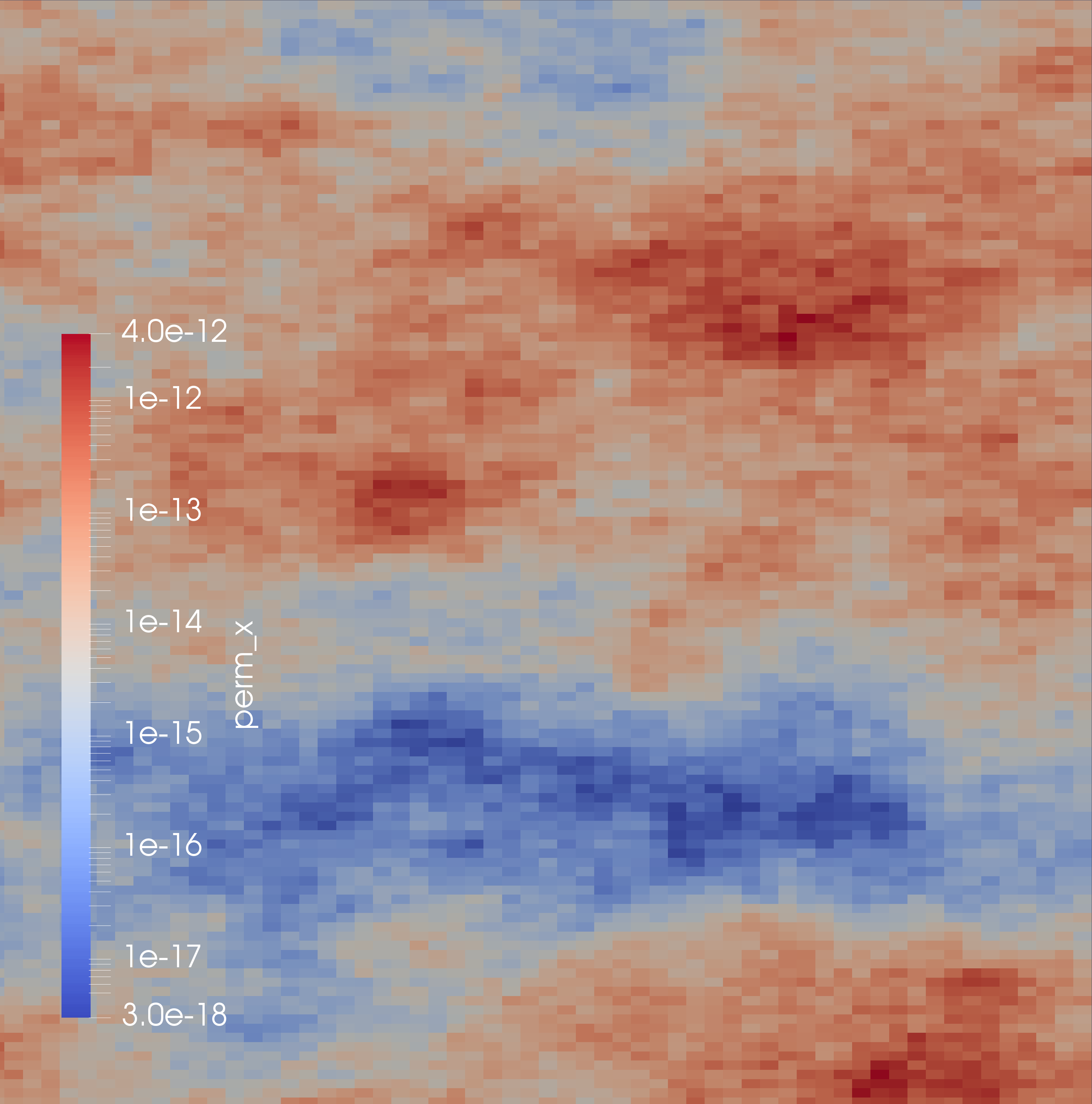}
 \caption{Log scale}
 \end{subfigure}
 \caption{Permeability of SPE10 test case (m$^2$).}\label{fig:perm_x}
\end{figure}

The domain is a square with dimensions 365.76$\times$365.76 meters, and the mesh is 60$\times$120. For permeability and porosity fields, we use the benchmark problem SPE10 \cite{christie2001tenth}. This problem has a highly heterogeneous permeability field. We consider a 60$\times$120 slice in the xy direction. The permeability, which is isotropic in the xy plane, is illustrated in Figure \ref{fig:perm_x}. We do not include gravity for the 2D simulations.

For the well case (W), we have one production well and one injection well. These are located in the upper half of the domain in the regions of high permeability. For the injection and production rates, we use the Peaceman well model. The bottom-hole pressure for the injection well is fixed at $6.895\times 10^7$ Pa, and $2.7579\times 10^7$ Pa for the production well. The maximum rate is set to $q=1.8\times 10^{-3}$m$^{3}$s$^{-1}$, although this is only achieved for the high permeability cases. The initial temperature in the reservoir is 288.706 K and the injection temperature is 422.039 K. For the heater case (H), heater placement is the same as for the well case, and so are the initial and heating temperatures. For the well and heater case (W+H), we combine both wells and heaters. For the high permeability cases (h.p.), we increase the permeability by a factor of 1,000. While the resulting permeability values are not representative of physical ones, they give a simple example of advection-dominated heat flow. 

For each case, we simulate injection and production for 1000 days where the time steps are chosen adaptively such that Newton's method converges in around 4 iterations. The average linear iterations per nonlinear iteration are shown in Table \ref{tab:spe10}.

\begin{table}[htb] \centering
\caption{SPE10 test cases. Average linear iterations per nonlinear iteration.}\label{tab:spe10}
 \begin{tabular}{|c c c c c c|}
 \hline
 method/case & W & H & W+H & h.p. W & h.p. W+H \\
 \hline
 Block & 5.88 & 5.42 & 6.60 & 14.5 & 14.0\\

 CPR   & 6.67 & 6.27 & 6.69 & 11.4 & 11.0\\
 \hline
 \end{tabular}
\end{table}

We observe that for the first three cases in Table \ref{tab:spe10}, the block preconditioner performs better than CPR in terms of the number of GMRES iterations, but that CPR performs best for the high permeability cases. The heat flow for the first three cases is diffusion-dominated, especially when the oil is not yet heated. For the high permeability cases, advection dominates. This change in performance appears later in the simulation when temperature has increased everywhere between the two wells. This indicates that CPR can still be a good choice if temperature is simply transported by the fluid flow. However, we will see in the next sections that the block preconditioner is a more scalable method.

\subsection{Numerical justification of the Schur complement approximation}

We now perform a numerical comparison of the action of the inverses of the different Schur complement approximations. We use the cases given in Section \ref{sec:spe10}. In Table \ref{tab:conds}, we compare the different Schur complement approximations by looking at the condition number of their inverse applied to the full Schur complement. While this condition number does not directly inform us about how well the preconditioner performs, it is a good indication of the quality of the approximations. For the cases, H and W stand for heaters and wells, respectively, and h.p. stands for high permeability (increased by a factor 1,000). We observe that $\tilde{S}_T$ is a good Schur complement approximation even for the high permeability cases where the other approximations struggle. 

\begin{table}[htb]
        \centering
\caption{Condition numbers (upper bounds) of the different matrices and Schur complement approximations for various cases}
\label{tab:conds}
\renewcommand*{\arraystretch}{1.15}
 \begin{tabular}{|c c c c c c|}
 \hline
 matrix/case & H & W & W+H & h.p. W & h.p. W+H \\
 \hline
 $\myrev{\tilde S_\mathrm{diag}^{-1}} S$   & 1.061  & 20.75 & 3.323   & 8.703e7& 2.191e7 \\
 $\myrev{\tilde S_\mathrm{A_{TT}}^{-1} }S$ & 1.063  & 28.08 & 4.277   & 4117   & 2467 \\
 ${\tilde S_T}^{-1} S$          & 1.023  & 1.097 & 1.1717  & 5.969  & 5.939 \\
 $A_{TT}$                   & 5.64e5 & 27.88 & 5.64e5  & 2324   & 2.143e5 \\
 $S$                        & 5.479e5& 2.862 & 5.717e5 & 20.60  & 4.976e5 \\
 \hline
 \end{tabular}
\end{table}

In terms of the performance of the solver, $\tilde{S}_T$ always results in fewer GMRES iterations (results not shown here). For harder cases (for example high permeability), this difference is significant; the linear solver can even fail to converge before the prescribed maximum number of iterations. \rev{In the next section, we will see that the other Schur complement approximations struggle in anisotropic medium.}

\subsection{Problem size scaling}\label{sec:probsizescaling}
We now investigate the performance of CPR and our block preconditioner as we refine a mesh. \rev{For two cases, we will also consider the Schur complement approximations $\tilde S_{A_{TT}}$ and $\tilde S_\mathrm{diag}$.} To this end, we test a case with homogeneous permeability and porosity fields. The domain is a square with dimensions 20$\times$ 20 meters and uniform porosity $\phi = 0.2$. We test both isotropic and anisotropic permeability fields. We refine the mesh from a $20\times20$ grid to $320\times 320$. 


We begin with an isotropic permeability of $3\times 10^{-13}$ m$^2$. For all cases, the injection/heating temperature is 422.039K. For all cases except Case III, the initial temperature is 288.706K. For each case, we take two time steps and calculate the average number of linear iterations per nonlinear iteration. For Case I-IV, the time step is 10 days, and for Case V, 12 hours.

For Case I, we have 6 heaters in the domain. In Table \ref{tab:isoH}, we observe that the number iteration increases by 9 times for CPR, while it increases by less than 50\% for the block preconditioner. 

\begin{table}[htb] \centering
\caption{Case I: Heater case in isotropic medium. Average linear iterations per nonlinear iteration.}\label{tab:isoH}
 \begin{tabular}{|c c c c c c|}
 \hline
 method/$N$ & 20 & 40 & 80 & 160 & 320 \\
 \hline
 Block & 2.57 & 3.23 & 2.86 & 3.44 & 3.71\\

 CPR   & 3.4 & 5.38 & 9.09 & 16.3 & 30.7\\
 
 \hline
 \end{tabular}
\end{table}

For Case II, we have injection wells and 3 production wells. The wells operate at constant injection and production rates of $q = 5 \times 10^{-8} m^3 s^{-1}$. In Table \ref{tab:isoW}, we observe that the number of iterations for CPR increases by 10 times while it only increases by less than 50 \% for the block preconditioner \rev{with the Schur complement approximation $\tilde S_T$. We observe a similar increase in iterations for the block preconditioner with the Schur complement approximations $\tilde{S}_{A_{TT}}$ and $\tilde{S}_\mathrm{diag}$.}

\begin{table}[htb] \centering
\caption{Case II: Well case in isotropic medium. Average linear iterations per nonlinear iteration.}\label{tab:isoW}
 \begin{tabular}{|c c c c c c|}
 \hline
 method/$N$ & 20 & 40 & 80 & 160 & 320 \\
 \hline
 Block & 2.43 & 2.43 & 2.86 & 3.28 & 3.71 \\

 CPR   & 3.71 & 5.71 & 9.86 & 19.4 & 37.4\\
 
  \rev{Block ($\tilde{S}_{A_{TT}}$)} &\rev{ 4.57 }&\rev{ 5    }&\rev{ 5.57 }&\rev{ 6.29 }&\rev{ 6.57} \\
 
 \rev{Block ($\tilde{S}_\mathrm{diag}$) }&\rev{ 4.14 }&\rev{ 4.43 }&\rev{ 5.29 }&\rev{ 5.86 }&\rev{ 6.43} \\
 \hline
 \end{tabular}
\end{table}

For Case III, we also have 3 injection wells and 3 production wells. To allow higher rates and faster flow, we increase the initial temperature to 320K. The wells operate at injection and production rates $q =  10^{-6} m^3 s^{-1}$. In Table \ref{tab:isoW}, we observe that the number of iterations for CPR increases by more than 10 times while it only increases by around 50 \% for the block preconditioner.
\begin{table}[htb] \centering
\caption{Case III: Higher injection well case in isotropic medium. Average linear iterations per nonlinear iteration.}\label{tab:isoHighW}
 \begin{tabular}{|c c c c c c|}
 \hline
 method/$N$ & 20 & 40 & 80 & 160 & 320 \\
 \hline
 Block & 3.67 & 4.38 & 4.7 & 5.10 & 5.52 \\

 CPR   & 4.71 & 7.31 & 13.1 & 24.7 & 50.6\\
 \hline
 \end{tabular}
\end{table}

For case IV and V, we increase the permeability in the $x$-direction to $3\times 10^{-11}$m$^2$. For Case IV, we have 6 heaters and observe the same trend as the previous case in Table \ref{tab:anisoH}.

\begin{table}[htb] \centering
\caption{Case IV: Heater case in anisotropic medium. Average linear iterations per nonlinear iteration.}\label{tab:anisoH}
 \begin{tabular}{|c c c c c c|}
 \hline
 method/$N$ & 20 & 40 & 80 & 160 & 320 \\
 \hline
 Block & 2.31 & 2.67 & 3.25 & 3.67 & 3.86\\

 CPR   & 3.11 & 4.56 & 8.56 & 15.8 & 30.4\\
 \hline
 \end{tabular}
\end{table}

For Case V, we have 3 injection wells and 3 production wells. The wells operate at constant injection and production rate $q = 1 \times 10^{-6} m^3 s^{-1}$. In this case, the flow is much faster and thus the time step size is reduced to half a day for the convergence of Newton's method. In Table \ref{tab:anisoW}, we observe that the number of iterations is doubled for the \myrev{b}lock preconditioner \rev{with $\tilde S_T$}, increased by 4 times for CPR\rev{, and slightly less for the block preconditioner with $\tilde S_{A_{TT}}$. Additionally, the block preconditioner with $\tilde S_\mathrm{diag}$ fails to converge within 200 GMRES iterations.}

\begin{table}[htb] \centering
\caption{Case \myrev{V}: Well case in anisotropic medium. Average linear iterations per nonlinear iteration.}\label{tab:anisoW}
 \begin{tabular}{|c c c c c c|}
 \hline
 method/$N$ & 20 & 40 & 80 & 160 & 320 \\
 \hline
 Block & 2.38 & 3.27 & 4.52 & 4.68 & 5.36 \\

 CPR   & 2.86 & 3.6 & 4.76 & 7.0 & 12.04\\
 
 \rev{Block ($\tilde{S}_{A_{TT}}$)}      & \rev{9} & \rev{17.1} & \rev{24.1} & \rev{27.9} & \rev{31.6} \\
 
 \rev{Block ($\tilde{S}_\mathrm{diag}$)} & \rev{$>200$} & \rev{$>200$} & \rev{$>200$}  &\rev{$>200$}  &\rev{$>200$} \\
 \hline
 \end{tabular}
\end{table}

In summary, as we refine the mesh, the number of iterations has a very small increase for the block preconditioner, but a large increase for CPR. The heat diffusion is much more noticeable on fine meshes, which CPR does not treat appropriately. However, coarser meshes are more common in commercial reservoir simulators. 

Note that the success of the block preconditioner is also due the linear scalability of AMG for elliptic problems. By removing the need for ILU, we get a nearly mesh-independent preconditioner.


\subsection{Parallel scaling}\label{sec:parallel}
We now compare the performance of the two methods in parallel. We look at both weak and strong scaling.

\subsubsection{Weak scaling}
For weak scaling, we compare the parallel performance of the methods as we increase the number of processors and problem size. The domain is $50\times 50\times 50$ meters with an $N\times N \times N$ grid. Since this is a 3D case, we include gravity. The permeability is $3\times10^{-13}$m$^2$ and the porosity is 0.2. We seek to have around \myrev{100,000} degrees of freedom per processor. Thus, for the number of processors 1, 2, 4, 8, and 16, we have $N$ = 36, 46, 58, 73, 92.

For the heating case, we have two sets of 21 heaters near the top and bottom of the domain. We take two time steps of 100 days and illustrate the results in Table \ref{tab:paralH}. We observe that the number of iterations increases by around 20 \% for the block preconditioner and triples for CPR. 

\begin{table}[htb] \centering
 \caption{Weak scaling: 3D Heating case. Average linear iterations per nonlinear iteration.}\label{tab:paralH}
 \begin{tabular}{|c c c c c c|}
 \hline
 method/num. proc. & 1 & 2 & 4 & 8 & 16 \\
 \hline
 Block & 7.5 & 7.9 & 8.25 & 8.75 & 9.29 \\

 CPR   & 15.75 & 22.3 & 29.9 & 38.5 & 45.6 \\
 \hline
 \end{tabular}
\end{table}

For the well case, we have 21 injection wells near the top of the domain, and 21 production wells near the bottom. All wells operate at a constant injection/production rate $q = 10^{-7} m^3 s^{-1}$. We take two time steps of 10 days and illustrate the performance of the methods in Table \ref{tab:paralW}. We observe that the number of iterations for the block preconditioner increases by around 40\% while the number of iterations for CPR nearly triples. 

\begin{table}[htb] \centering
 \caption{Weak scaling: 3D Well case. Average linear iterations per nonlinear iteration.}\label{tab:paralW}
 \begin{tabular}{|c c c c c c|}
 \hline
 method/num. proc. & 1 & 2 & 4 & 8 & 16 \\
 \hline
 Block & 4.29 & 4.43 & 4.71 & 5.25 & 5.89 \\

 CPR   & 6.57 & 10.0 & 12.3 & 16.3 & 18.4 \\
 \hline
 \end{tabular}
\end{table}

\subsubsection{Strong scaling}
We use the same problem as the previous section on the finest mesh. We keep the problem size fixed while increasing the number of processors. For reservoir simulation, strong scaling is usually more relevant than weak scaling. Indeed, reservoir models often come with a (usually rather coarse) fixed grid. As observed in Section \ref{sec:probsizescaling}, CPR does not behave as well on a fine mesh. Thus, keeping the mesh size constant is a good way of isolating the parallel performance of the methods. 

In Tables \ref{tab:strongH} and \ref{tab:strongW}, we illustrate the strong scaling results for the heating and well cases, respectively. For the block preconditioner, we observe that the number of iteration is essentially independent of the number of processors used. This is thanks to the parallel capability of BoomerAMG. On the other hand, the number of iterations for CPR exhibit a small but progressive increase. This is because the second stage of CPR uses Block ILU, which becomes a weaker preconditioner as the number of blocks increases. Therefore, this trend will continue as the number of processors increases. 

\begin{table}[htb] \centering
 \caption{Strong scaling: 3D Heating case. Average linear iterations per nonlinear iteration.}\label{tab:strongH}
 \begin{tabular}{|c c c c c c|}
 \hline
 method/num. proc. & 1 & 2 & 4 & 8 & 16 \\
 \hline
 Block & 8.75 & 8.57 & 8.57 & 9.43 & 9.29 \\

 CPR   & 38.0 & 43.3 & 44.0 & 44.7 & 45.6 \\
 \hline
 \end{tabular}
\end{table}

\begin{table}[htb] \centering
 \caption{Strong scaling: 3D Well case. Average linear iterations per nonlinear iteration.}\label{tab:strongW}
 \begin{tabular}{|c c c c c c|}
 \hline
 method/num. proc. & 1 & 2 & 4 & 8 & 16 \\
 \hline
 Block & 6.22 & 5.44 & 5.78 & 6.11 & 5.89 \\

 CPR   & 14.0 & 16.1 & 16.7 & 17.9 & 18.4 \\
 \hline
 \end{tabular}
\end{table}

\section{Conclusion}\label{sec:conclusion}

In this work, we have implemented a fully implicit parallel non-isothermal porous media flow simulator including two preconditioning strategies, CPR and \myrev{a block preconditioner with our own Schur complement approximation}. We have tested the performance of these methods as preconditioners for GMRES. \rev{Our Schur complement approximation performs better than simple one, especially in cases with heterogeneous or anisotropic permeability.} While the block preconditioner performs well for diffusion-dominated cases, CPR is still the method of choice for advection-dominated (manufactured) cases, at least in serial. However, the block preconditioner scales optimally with problem size while CPR does not do well under mesh refinement. Additionally, the block preconditioner remains efficient in parallel, while the CPR iteration count increases gradually as we increase the number of processors. 

This research demonstrates that a preconditioning strategy which considers the diffusive effect of temperature is important for diffusion-dominated cases. In non-isothermal multiphase flow, the energy equation is treated in CPR like a hyperbolic equation. A coupled solution of pressure and temperature using multigrid is key to \rev{methods for multiphase flow currently being developed.}


\section{Acknowledgments}

This publication is based on work partially supported by the EPSRC Centre For Doctoral Training in Industrially Focused Mathematical Modelling (EP/L015803/1) in collaboration with Schlumberger.

\appendix

 \bibliographystyle{elsarticle-num} 
 \bibliography{bibsinglethermal}





\end{document}